\pdfoutput=1
\documentclass[11pt,a4paper]{article}

\usepackage[latin1]{inputenc}
\usepackage{eucal,amssymb,amsfonts,amsmath,amsthm,graphicx,bbm,epsfig,color,fullpage,comment,bbold}
\newtheorem{Theorem}{Theorem}
\newtheorem{Definition}[Theorem]{Definition}

\newtheorem{Lemma}[Theorem]{Lemma} 
\newtheorem{Proposition}[Theorem]{Proposition}
\newtheorem{Conjecture}[Theorem]{Conjecture}
 \newtheorem{Open
  question}[Theorem]{Open question}

\begin{document}

\title{Confidence intervals for the critical value in the divide and color model}

\author{Andr{\'a}s B{\'a}lint\thanks{Chalmers University of Technology,
    e-mail: {abalint@few.vu.nl}} \and
  Vincent Beffara\thanks{UMPA-ENS Lyon,
    e-mail: {vbeffara@ens-lyon.fr}} \and
  Vincent Tassion\thanks{UMPA-ENS Lyon,
    e-mail: {vincent.tassion@ens-lyon.fr}}}

\maketitle

\begin{abstract}
  We obtain confidence intervals for the location of the percolation
  phase transition in H\"aggstr\"om's divide and color model on the
  square lattice $\mathbb{Z}^2$ and the hexagonal lattice
  $\mathbb{H}$. The resulting probabilistic bounds are much tighter
  than the best deterministic bounds up to date; they give a clear
  picture of the behavior of the DaC models on $\mathbb{Z}^2$ and
  $\mathbb{H}$ and enable a comparison with the triangular lattice
  $\mathbb{T}$. In particular, our numerical results suggest
  similarities between DaC model on these three lattices that are in
  line with universality considerations, but with a remarkable
  difference: while the critical value function $r_c(p)$ is known to
  be constant in the parameter $p$ for $p<p_c$ on $\mathbb{T}$ and
  appears to be linear on $\mathbb{Z}^2$, it is almost certainly
  non-linear on $\mathbb{H}$.
\end{abstract}



  \section{Introduction}

  Our object of study in this paper is the critical value function in
  H\"aggstr\"om's divide and color (DaC) model \cite{HaggstromDaC}.
  This is a stochastic model that was originally motivated by physical
  considerations (see \cite{HaggstromDaC,CLM}), but it has since then
  been used for biological modelling in \cite{GPG} as well and
  inspired several generalizations (see, \emph{e.g.},
  \cite{HH,BCM,GrGr}).  Our results concerning the location of the
  phase transition give a clear picture of the behavior of the DaC
  model on two important lattices and lead to intriguing open
  questions.

  Our analysis will be based on the same principles as \cite{RW},
  where confidence intervals were obtained for the critical value of
  Bernoulli bond and site percolation on the 11 Archimedean lattices
  by a modification of the approach of \cite{BBW}.  The main idea in
  \cite{BBW,RW} is truly multidisciplinary and attractive, namely to
  reduce a problem which has its roots in theoretical physics by deep
  mathematical theorems to a situation in which a form of statistical
  testing by numerical methods becomes possible.  Our other main goal
  with this paper is to demonstrate the strength of this strategy by
  applying it to a system which is essentially different from those in
  its previous applications. In particular, in the DaC model, as
  opposed to the short-range dependencies in \cite{BBW} and the
  i.i.d.\ situation in \cite{RW}, one has to deal with correlations
  between sites at arbitrary distances from each other. We believe
  that the method of \cite{BBW,RW} has a high potential to be used in
  a number of further models (see \emph{e.g.}~\cite{DHP} where a very
  similar approach is followed) and deserves higher publicity than it
  enjoys at the moment.

  Given a graph $G$ with vertex set $\mathcal{V}$ and edge set
  $\mathcal{E}$ and parameters $p,r\in [0,1]$, the DaC model on $G$ is
  defined in two steps: first, Bernoulli bond percolation with density
  $p$ is performed on $G$, and then the resulting open clusters are
  independently colored black (with probability $r$) or white (a more
  detailed definition will follow in the next paragraph).  Note that
  this definition resembles the so-called random-cluster (or FK)
  representation of the ferromagnetic Ising model, with two important
  differences: a product measure is used in the DaC model in the first
  step instead of a random-cluster measure with cluster weight $2$ and
  the second step is more general here in that all $r\in [0,1]$ are
  considered instead of only $1/2$.

  Now we set the terminology that is used throughout, starting with an
  alternative (equivalent) definition of the DaC model which goes as
  follows.  First, an \emph{edge configuration} $\eta \in
  \{0,1\}^{\mathcal{E}}$ is drawn according to the product measure
  $\nu _p^{\mathcal{E}}$ where $\nu _p$ is the probability measure on
  $\{0,1\}$ with $\nu _p(\{ 1\})=1-\nu _p(\{ 0\})=p$.  In the second
  step, a \emph{site configuration} $\xi \in \{0,1\}^{\mathcal{V}}$ is
  chosen by independently assigning state $1$ with probability $r$ or
  otherwise $0$ to each vertex, conditioning on the event that there
  exists no edge $e=\left< v,w\right> \in \mathcal{E}$ such that $\eta
  (e)=1$ and $\xi (v)\neq \xi (w)$. We denote the probability measure
  on $\{0,1\}^{\mathcal{V}}\times \{0,1\}^{\mathcal{E}}$ associated to
  this procedure by $\mathbb{P}_{p,r}^G$. An edge $e$ (a vertex $v$)
  is said to be \emph{open} or \emph{closed} (\emph{black} or
  \emph{white}) if and only if it is in state $1$ or $0$,
  respectively. We will call the maximal subsets of $\mathcal{V}$
  connected by open edges \emph{bond clusters}, and the maximal
  monochromatic connected (via the edge set of $\mathcal{E}$, not only
  the open edges!) subsets of $\mathcal{V}$ \emph{black} or
  \emph{white clusters}. We write $C_v(\eta )$ for the bond cluster of
  a vertex $v$ in the edge configuration $\eta $ and use $\Omega _S$
  to denote $\{0,1\}^S$ for arbitrary sets $S$.

  Note that the measure $\mathbb{P}_{p,r}^G$ is concentrated on the
  set of pairs $(\eta, \xi)$ such that for all edges $e=\left<
    v,w\right> \in \mathcal{E}$, $\xi (v)=\xi(w)$ whenever
  $\eta(e)=1$. When this compatibility condition is satisfied, we
  write $\eta \sim \xi$.

  \medskip

  For infinite graphs $G$, there are two types of phase transitions
  present in the DaC model in terms of the appearance of infinite
  $1$-clusters; first, there exists $p_c=p_c^G\in [0,1]$ such that
  $\mathbb{P}_{p,r}^G($there exists an infinite bond cluster$)$ is $0$
  for $p<p_c$ and $1$ for $p>p_c$. Second, for each fixed $p$, there
  exists $r_c=r_c^G(p)$ such that $\mathbb{P}_{p,r}^G($there exists an
  infinite black cluster$)$ is $0$ for $r<r_c$ and \emph{positive} for
  $r>r_c$. For more on the different character of these two types of
  phase transitions, see \cite{BBT}.  A key feature of the DaC model
  (as noted in \cite{HaggstromDaC}) is that while it is close in
  spirit to the Ising model, its simulation is straightforward from
  the definition and does not require sophisticated MCMC
  algorithms. In this paper, we will exploit this feature in order to
  learn about the values and various features of the \emph{critical
    value function} $r_c^G(p)$.

\begin{figure}[ht]
    \centering
    \includegraphics[scale=0.35, trim= 0mm 0mm 0mm 0mm, clip
    ]{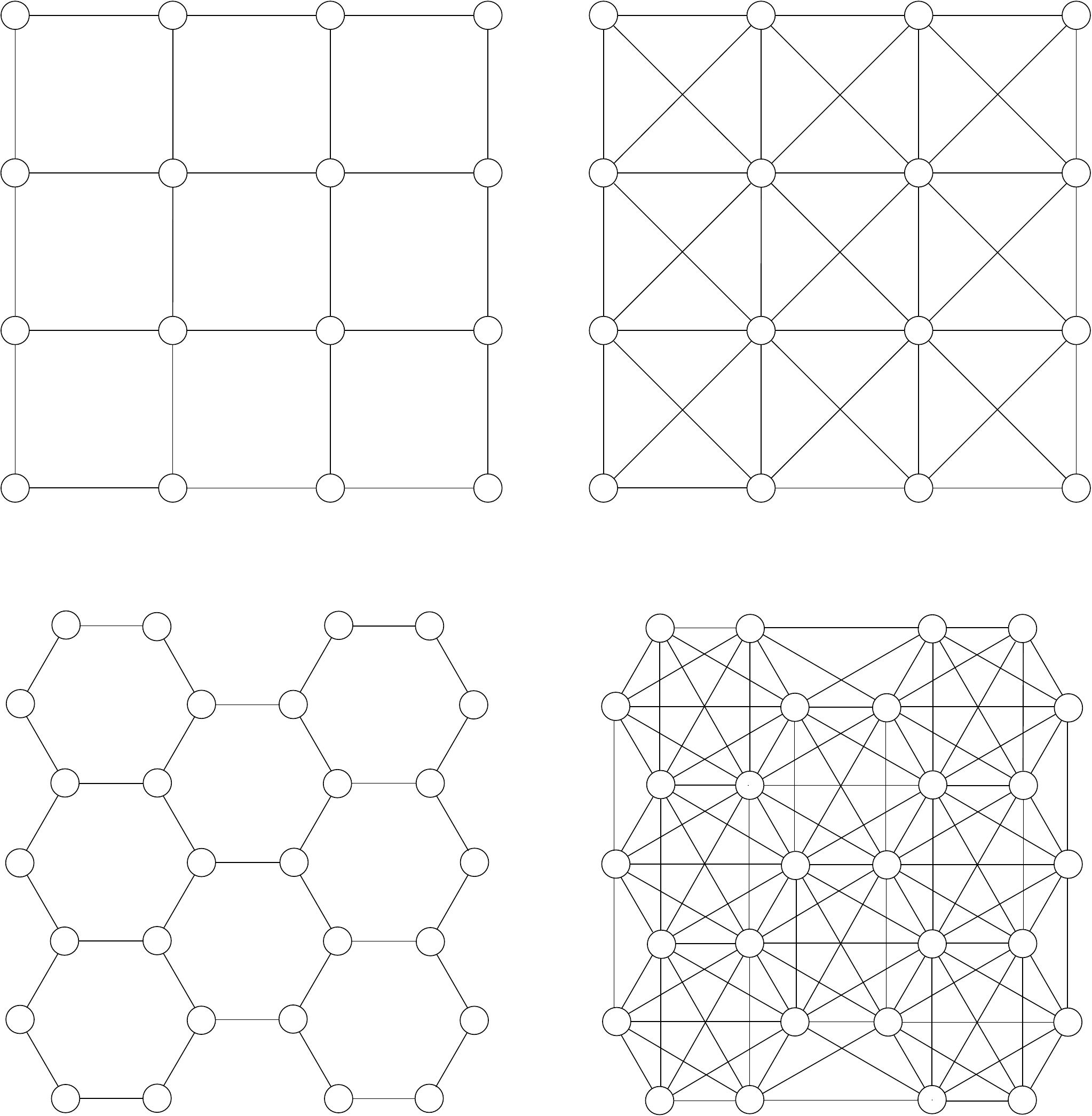}
    \caption{A finite sublattice of the square lattice $\mathbb{Z}^2$
      (above left) and the hexagonal lattice $\mathbb{H}$ (below left)
      and their respective matching lattices (right).}
    \label{lattices}
  \end{figure}

  Monotonicity and continuity properties of the function $r_c^G(p)$
  for general graphs have been studied in \cite{BBT}. Here we will
  focus on two specific graphs, namely the \emph{square lattice}
  $\mathbb{Z}^2$ and the \emph{hexagonal lattice} $\mathbb{H}$ (see
  Figure \ref{lattices}), for which $p_c^{\mathbb{Z}^2}=1/2$ and
  $p_c^{\mathbb{H}}=1-2\sin (\pi /18)\approx 0.6527$ (see
  \cite{kesten-book}). Our reason for this restriction is twofold:
  first, these two are the most commonly considered planar lattices
  (apart from the triangular lattice $\mathbb{T}$, for which the
  critical value function $r_c^{\mathbb{T}}$ has been completely
  characterized in \cite{BCM}), whence results about these cases are
  of the greatest interest.  On the other hand, the DaC model on these
  lattices enjoys a form of duality (described in Section
  \ref{duality-subsection}) which is a key ingredient for the analysis
  we perform in this paper.

  Fixing $\mathcal{L}\in \{\mathbb{Z}^2,\mathbb{H}\}$, it is trivial
  that $r_c^{\mathcal{L}}(p)=0$ for all $p>p_c^{\mathcal{L}}$, and it
  easily follows from classical results on Bernoulli bond percolation
  that $r_c^{\mathcal{L}}(p_c^{\mathcal{L}})=1$ (see \cite{BCM} for
  the case $\mathcal{L}=\mathbb{Z}^2$). However, there are only very
  loose theoretical bounds for the critical value when
  $p<p_c^{\mathcal{L}}$: the duality relation (\ref{rcplusrcstaris1})
  in Section \ref{duality-subsection} below and renormalization
  arguments as in the proof of Theorem 2.6 in \cite{HaggstromDaC} give
  that $1/2\leq r_c^{\mathcal{L}}(p)<1$ for all such $p$, and
  Proposition 1 in \cite{BBT} gives just a slight improvement of these
  bounds for very small values of $p$.  Therefore, our ultimate goal
  here is to get good estimates for $r_c^{\mathcal{L}}(p)$ with
  $p<p_c^{\mathcal{L}}$.

  We end this section with an outline of the paper. Section
  \ref{fsc-section} contains a crucial reduction of the
  infinite-volume models to a \emph{finite} situation by a criterion
  that is stated in terms of a finite sublattice but nonetheless
  implies the existence of an infinite cluster. This method, often
  called static renormalisation in percolation, is a particular
  instance of coarse graining.  We then describe in Section
  \ref{confidence-interval-section} how the occurrence of this finite
  size criterion can be tested in an efficient way and obtain
  confidence intervals for $r_c^{\mathcal{L}}(p)$ as functions of
  uniform random variables (Proposition \ref{confintprop}). Finally,
  we implement this method using a (pseudo)random number generator,
  and present and discuss the numerical results in Section
  \ref{results-section}.

  \section{Finite size criteria}
  \label{fsc-section}

  \subsection{An upper bound for $r_c(p)$}
  \label{upperbound-subsection}

  In this section, we will show how to obtain an upper bound for
  $r_c^{\mathcal{L}}(p)$ by deducing a finite size criterion for
  percolation in the DaC model (Proposition
  \ref{fscpropositionforp>0}).  This criterion, which is a
  quantitative form of Lemma 2.10 in \cite{BCM}, will play a key role
  in Sections
  \ref{confidence-interval-section}--\ref{results-section}.  To
  enhance readability, we will henceforth focus on the case
  $\mathcal{L}=\mathbb{Z}^2$ and mention $\mathcal{L}=\mathbb{H}$ only
  when the analogy is not straightforward. Accordingly, we will write
  $\mathbb{P}_{p,r}$ and $r_c(p)$ for
  $\mathbb{P}^{\mathbb{Z}^2}_{p,r}$ and $r_c^{\mathbb{Z}^2}(p)$
  respectively, and denote the edge set of $\mathbb{Z}^2$ by
  $\mathcal{E}^2$.  Let us first recall a classical result (Lemma
  \ref{1deplemma} below) concerning 1-dependent percolation.

\begin{Definition}
  Given a graph $G=(\mathcal{V},\mathcal{E})$, a probability measure
  $\nu $ on $\{0,1\}^{\mathcal{E}}$ is called \textit{1-dependent} if,
  whenever $S \subset \mathcal{E}$ and $T \subset \mathcal{E}$ are
  vertex-disjoint edge sets, the state of edges in $S$ is independent
  of that of edges in $T$ under $\nu $.
\end{Definition}

It follows from standard arguments or from a general theorem of
Liggett, Schonmann and Stacey \cite{LSS} that if each edge is open
with a sufficiently high probability in a 1-dependent bond percolation
on $\mathbb{Z}^2$, then the origin is with positive probability in an
infinite bond cluster. Currently the best bound is given by Balister,
Bollob\'as and Walters \cite{BBW}:

\begin{Lemma}\emph{(\cite{BBW})}\label{1deplemma}
  Let $\nu $ be any 1-dependent bond percolation measure on
  $\mathbb{Z}^{2}$ in which each edge is open with probability at
  least $0.8639$. Then the probability under $\nu $ that the origin
  lies in an infinite bond cluster is positive.
\end{Lemma}

\begin{figure}[ht]
  \begin{center}
    \includegraphics{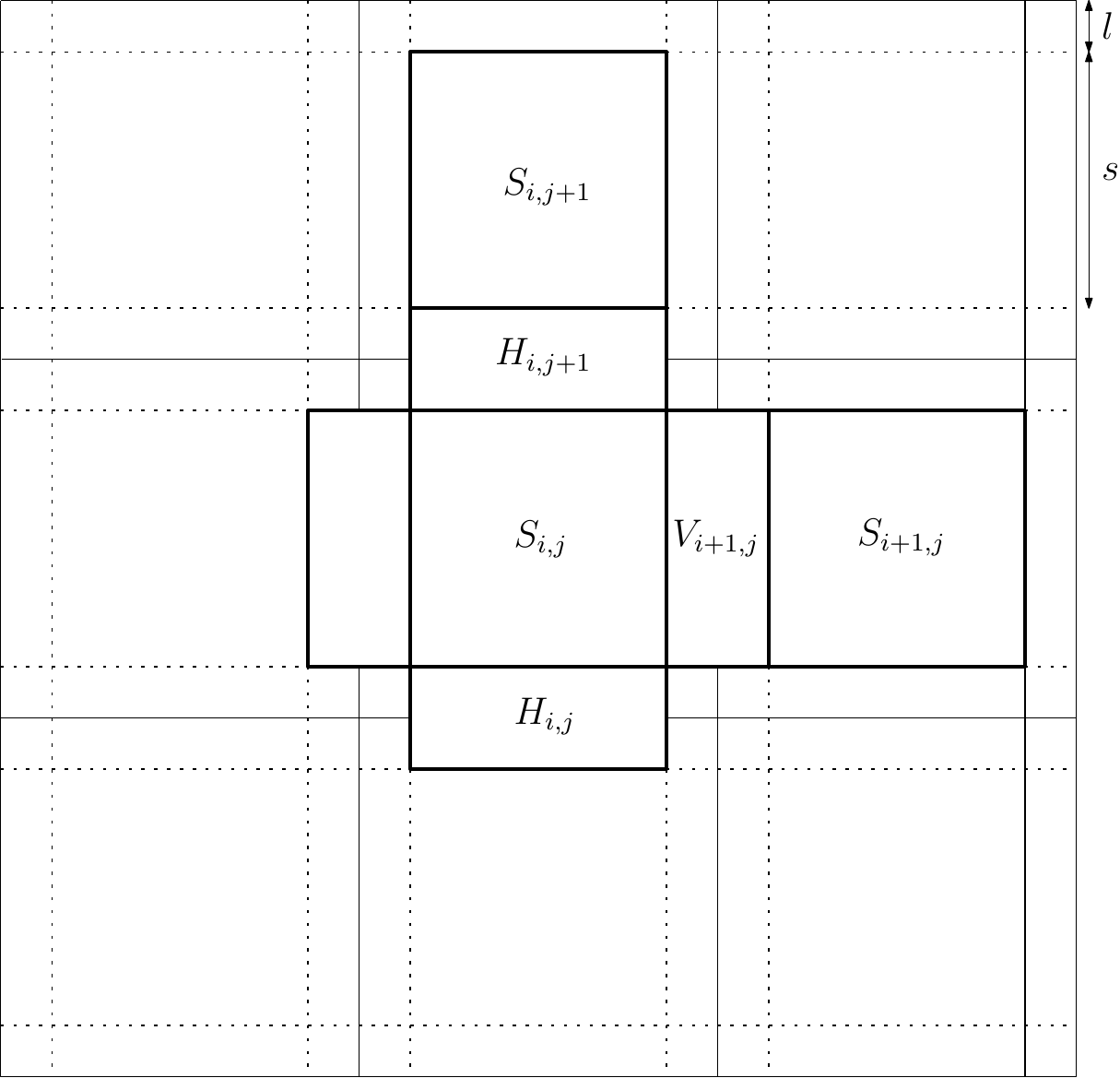}
    \caption{A partition of $\mathbb{R}^2$.}
    \label{partition}
  \end{center}
\end{figure}

Now, suppose that the lattice $\mathbb{Z}^2$ is embedded in the plane
the natural way (so that $v=(i,j)\in \mathbb{Z}^2$ has coordinates $i$
and $j$).  We consider the following partition of $\mathbb{R}^2$ (see
Figure \ref{partition}): given parameters $s \in \mathbb{N}=\{
1,2,\ldots \}$ and $\ell \in \mathbb{N}$, we take $k=s+2\ell $ and
define, for all $i,j\in \mathbb{Z}$, the $s\times s$ squares
$$
S_{i,j}=[ik+\ell ,ik+\ell +s]\times [jk+\ell ,jk+\ell +s],
$$
the $s \times 2\ell $ rectangles
$$
H_{i,j}=[ik+\ell ,ik+\ell +s]\times [jk-\ell ,jk+\ell ],
$$
the $2\ell \times s$ rectangles
$$
V_{i,j}=[ik-\ell ,ik+\ell ]\times [jk+\ell ,jk+\ell +s],
$$
and what remains are the $2\ell \times 2\ell $ squares $ [ik-\ell
,ik+\ell ]\times [jk-\ell ,jk+\ell ].  $

\begin{figure}[ht]
  \begin{center}
    \includegraphics[width=0.8\hsize]{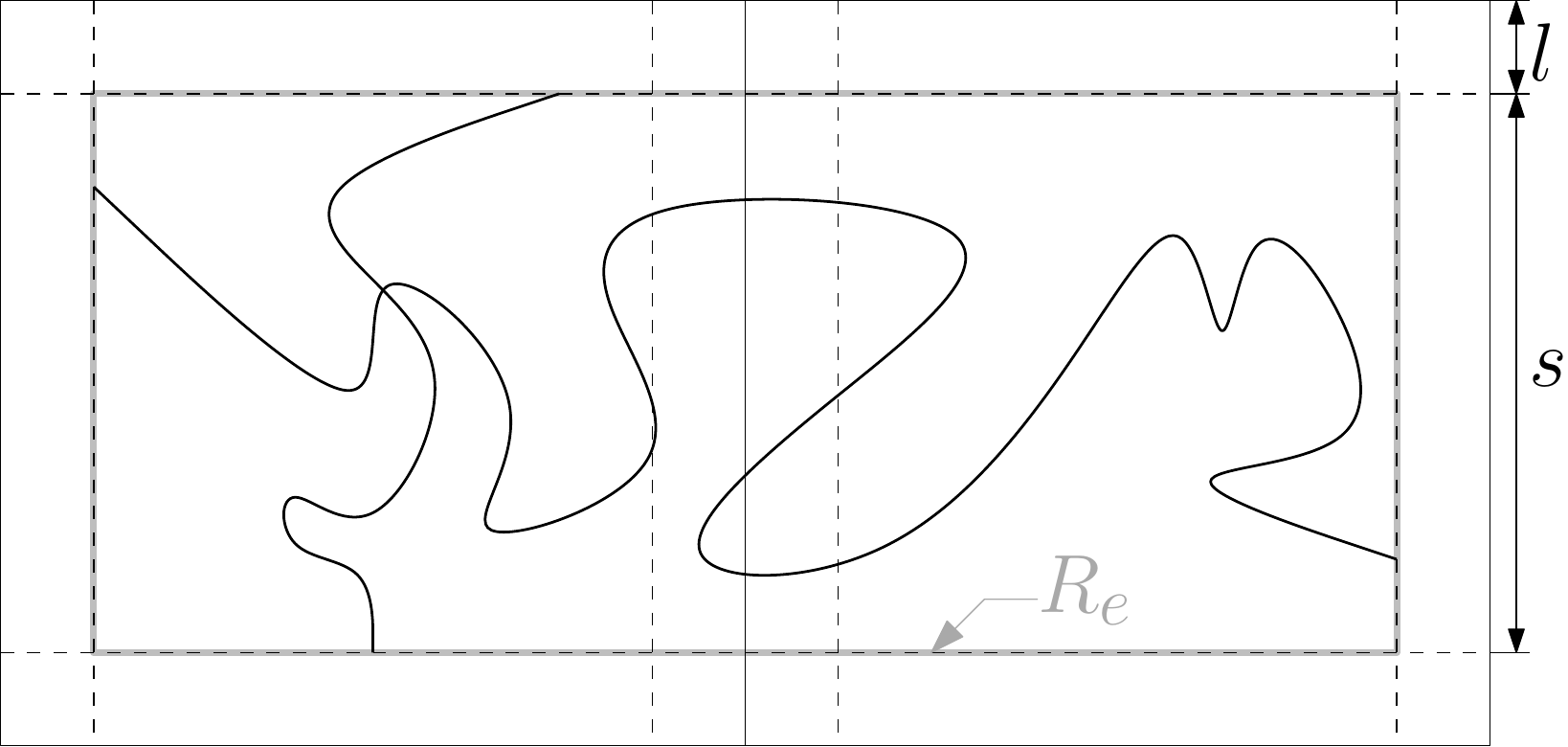}
    \caption{A black component in $R_e$ witnesses the occurrence of
      $E_e$.}
    \label{blackpaths}
  \end{center}
\end{figure}

We will couple $\mathbb{P} _{p,r}$ to a $1$-dependent bond percolation
measure. Define $f:\Omega _{\mathcal{E}^2}\times \Omega _{
  \mathbb{Z}^{2} }\rightarrow \Omega _{\mathcal{E}^2}$, as follows.
To each horizontal edge $e=\left< (i,j),(i+1,j)\right> \in
\mathcal{E}^2$, we associate a $(2\ell +2s) \times s$ rectangle
$R_e=S_{i,j}\cup V_{i+1,j}\cup S_{i+1,j}$ and the event $E_e$ that
there exists a left-right black crossing in $R_e$ (\emph{i.e.}, a
connected path of vertices all of which are black which links the left
side of $R_e$ to its right side) and an up-down black crossing in
$S_{i,j}$ (see Figure \ref{blackpaths}).
Here and below, a vertex in the corner of a rectangle is understood to
link the corresponding sides in itself.  For each vertical edge
$e=\left< (i,j),(i,j+1)\right> \in \mathcal{E}^2$, we define the $s
\times (2\ell +2s)$ rectangle $R_e=S_{i,j}\cup H_{i,j+1}\cup
S_{i,j+1}$ and the event $E_e=\{$up-down black crossing in $R_e$ and
left-right black crossing in $S_{i,j}\} \subset \Omega
_{\mathcal{E}^2}\times \Omega _{ \mathbb{Z}^{2} }$.  For each edge
$e\in \mathcal{E}^2$, we also consider the event $F_{e}=\{$there
exists a bond cluster which contains a vertex in $R_e$ and a vertex at
graph distance at least $\ell$ from $R_e\}\subset \Omega
_{\mathcal{E}^2}\times \Omega _{ \mathbb{Z}^{2} }$, and define $\tilde
{E}_e=E_e\cap F_e^c$.  Now for each configuration $\omega =(\eta ,\xi
) \in \Omega _{\mathcal{E}^2}\times \Omega _{ \mathbb{Z}^{2} }$, we
determine a corresponding bond configuration $f(\omega )=\gamma \in
\Omega _{\mathcal{E}^2}$ as follows: for all $e\in \mathcal{E}^2$, we
declare $e$ open if and only if $\tilde {E}_{e}$ holds (\emph{i.e.},
we define $\gamma (e) = 1$ if and only if $\omega \in \tilde{E}_{e}
$).  Finally, we define the probability measure $\nu = f _\ast
\mathbb{P} _{p,r}$ on $\Omega_{\mathcal{E}^2}$.

It is not difficult to check that $\nu $ is a 1-dependent bond
percolation measure.  Indeed, if $e$ and $e^{\prime }$ are two
vertex-disjoint edges in $\mathcal{E}^2$, then the corresponding
rectangles $R_{e}$ and $R_{e^{\prime }}$ are at graph distance at
least $2\ell $ from one another, hence $F_e^c$ and $F_{e^{\prime }}^c$
are independent.  Given that $F_e$ and $F_{e^{\prime }}$ do not hold,
the bond clusters in $R_e$ and $R_{e^{\prime }}$ are colored
independently of each other.  Keeping this in mind, a short
computation proves the independence of $\tilde{E}_e$ and
$\tilde{E}_{e^{\prime }}$ under $\mathbb{P}_{p,r}$, which implies the
1-dependence of $\nu $.

Note also that the function $f$ was chosen in such a way that if
$\gamma =f(\omega )\in \Omega _{\mathcal{E}^2}$ contains an infinite
open bond cluster, then $\omega$ contains an infinite black
cluster. Such configurations have zero $\mathbb{P} _{p,r}$-measure for
$r<r_c(p)$.  Finally, note that $\mathbb{P} _{p,r}(\tilde{E}_e)$ is
the same for all edges $e\in \mathcal{E}^2$. These observations
combined with Lemma \ref{1deplemma} imply that, denoting $\left<
  (0,0),(0,1)\right> \in \mathcal{E}^2$ by $e_1$, we have the
following result.

\begin{Proposition}\label{fscpropositionforp>0}
  Given any values of the parameters $s,\ell \in \mathbb{N}$, if $p$
  and $r$ are such that
  \begin{equation}\label{fscforp>0}
    \mathbb{P}_{p,r}( \tilde{E}_{e_1} ) \geq 0.8639,
  \end{equation}
  then $r_{c}(p)\leq r$.
\end{Proposition}

Note that Proposition \ref{fscpropositionforp>0} is indeed a finite
size criterion since the event $\tilde {E}_{e_1}$ depends on the state
of a finite number of edges and the color of a finite number of
vertices.  A similar criterion, which will imply a lower bound for
$r_c(p)$, will be given in Section \ref{lowerbound-subsection}.

\subsection{Duality}\label{duality-subsection}

A concept that is essential in understanding site percolation models
on $\mathcal{L}\in \{ \mathbb{Z}^2,\mathbb{H}\}$ is that of the
\emph{matching lattice} $\mathcal{L}^*$ which is a graph with the same
vertex set, $\mathcal{V}$, as $\mathcal{L}$ but more edges: the edge
set $\mathcal{E}^*$ of $\mathcal{L}^*$ consists of all the edges in
$\mathcal{E}$ plus the diagonals of all the faces of $\mathcal{L}$
(see Figure~\ref{lattices}).  The finiteness of a monochromatic
cluster in $\mathcal{L}$ can be rephrased in terms of circuits of the
opposite color in $\mathcal{L}^*$ and vice versa; see
\cite{kesten-book} for further details.  We say that $B\subset
\mathcal{V}$ is a \emph{black} $*$\emph{-component} in a color
configuration $\xi \in \Omega _{\mathcal{V}}$ if it is a black
component in terms of the lattice $\mathcal{L}^*$ (\emph{i.e.},
$\xi(v)=1$ for all $v\in B$ and $B$ is connected via $\mathcal{E}^*$).

Accordingly, there is yet another phase transition in the DaC model on
$\mathcal{L}$ at the point where an infinite black $*$-component
appears; formally, for each fixed $p\in [0,1]$, one can define
$r_c^*(p,\mathcal{L})$ as the value such that
$\mathbb{P}^{\mathcal{L}}_{p,r}($there exists an infinite black
$*$-component$)$ is $0$ for $r<r_c^*(p,\mathcal{L})$ and positive for
$r>r_c^*(p,\mathcal{L})$.  It was proved in \cite{BCM} that there is
an intimate connection between all the critical values in the DaC
model that we mentioned so far; namely, for all $p<p_c^{\mathcal{L}}$,
\begin{equation}\label{rcplusrcstaris1}
  r_c^{\mathcal{L}}(p)+r_c^*(p,\mathcal{L})=1.
\end{equation}
Actually, this relation was proved only for
$\mathcal{L}=\mathbb{Z}^2$, but essentially the same proof gives the
result for $\mathcal{L}=\mathbb{H}$ as well. The importance of this
result here is that due to the duality relation
(\ref{rcplusrcstaris1}), a lower bound for $r_c^{\mathcal{L}}(p)$ may
be obtained by giving an upper bound for $r_c^*(p,\mathcal{L})$.

\subsection{A lower bound for $r_c(p)$}\label{lowerbound-subsection}

As in Section \ref{upperbound-subsection}, we will focus on
$\mathcal{L}=\mathbb{Z}^2$ since the case $\mathcal{L}=\mathbb{H}$ is
analogous; we denote $r_c^*(p,\mathbb{Z}^2)$ here and in the next
section by $r_c^*(p)$.  Obviously $r_c(p)$ itself is an upper bound
for $r_c^*(p)$.  However, a better bound may be obtained by a slight
modification of the approach given in Section
\ref{upperbound-subsection}. For each $e\in \mathcal{E}^2$, let $R_e$
and $F_e$ be as in Section \ref{upperbound-subsection}, define $E_e^*$
by substituting black $*$-component for black component in the
definition of $E_e$, and take $\tilde {E}_e^*=E_e^*\cap F_e^c$. Then,
by similar arguments as those before Proposition
\ref{fscpropositionforp>0} and using (\ref{rcplusrcstaris1}), we get
the following:

\begin{Proposition}\label{fscstarpropositionforp>0}
  Given any values of the parameters $s,\ell \in \mathbb{N}$, if $p$
  and $r$ are such that
  \begin{equation}\label{fscstarforp>0}
    \mathbb{P}_{p,r}( \tilde{E}^*_{e_1} ) \geq 0.8639,
  \end{equation}
  then $r_{c}^*(p)\leq r$, and hence $r_c(p)\geq 1-r$.
\end{Proposition}

\section{The confidence interval} \label{confidence-interval-section}

The main idea in \cite{BBW,RW} is to reduce a stochastic model to a
new model in finite volume by criteria similar in spirit to those in
Section \ref{fsc-section} and do repeated (computer) simulations of
the new model to test whether the corresponding criteria hold.  The
point is that after a sufficiently large number of simulations, one
can see with an arbitrarily high level of confidence whether or not
the probability of an event exceeds a certain threshold.  By the
special nature of the events in question, statistical inferences
regarding the original, infinite-volume model may be made from the
simulation results.

To be able to follow this strategy, we will have to refine
Propositions
\ref{fscpropositionforp>0}--\ref{fscstarpropositionforp>0} as those
are concerned with the state of finitely many objects, but still in
the infinite-volume model. The adjusted criteria that truly are of
finite size are given below, see (\ref{GtildeEe}) and
(\ref{GtildeEestar}). Finding an efficient way of performing the
simulation step involves further obstacles. The main problem is that
it would be unfeasible to run a large number of separate simulations
for different values of $r$ to find, for a fixed $p$, the lowest value
of $r$ such that both (\ref{GtildeEe}) and (\ref{GtildeEestar}) seem
sufficiently likely to hold. We will tackle this difficulty with a
stochastic coupling, which is the simultaneous construction of several
stochastic models on the same probability space. Such a construction
will enable us to deal with \emph{all} values of $r\in [0,1]$ at the
same time and is very related to the model of invasion percolation.

After the description of the coupling, a ``theoretical'' confidence
interval (meaning a confidence interval as a function of i.i.d.\
random variables) for $r_c(p)$ is given in Proposition
\ref{confintprop}. The numerical confidence intervals obtained by this
method using computer simulations will be presented in Section
\ref{results-section}.  Note also that the inequalities
(\ref{GtildeEe}) and (\ref{GtildeEestar}) implicitly involve the
parameters $s$ and $\ell $ whose choices may influence the width of
the confidence intervals obtained; this issue is addressed before the
proof of Proposition \ref{confintprop}.  Our methods in this section
work for a general $p\in [0,p_c^{\mathcal{L}})$; we note that
substantial simplifications are possible in the case $p=0$
(\emph{i.e.}, in the absence of correlations), see \cite{RW}.

Fix $p\in [0,1/2)$ and $s, \ell \in \mathbb{N}$, and define the
rectangle $\tilde R_{e_1}=[0,2s+4\ell ]\times [0,s+2\ell ]$. Note that
for a configuration $\omega \in \Omega _{\mathcal{E}^2}\times \Omega
_{\mathbb{Z}^2}$, one can decide whether $\omega \in \tilde{E}_{e_1} $
(respectively $\omega \in \tilde{E}^*_{e_1} $) holds by checking the
restriction of $\omega $ to $\tilde R_{e_1}$. In fact, defining
$\tilde {G}=(\tilde {\mathcal{V}},\tilde {\mathcal{E}})$ as the
minimal subgraph of $\mathbb{Z}^2$ which contains $\tilde R_{e_1}$ and
considering the DaC model on $\tilde {G}$, it is easy to see that for
any $r\in [0,1]$, \smash{$\mathbb{P}_{p,r}^{\tilde {G}}(
  \tilde{E}_{e_1} )=\mathbb{P}_{p,r}( \tilde{E}_{e_1} )$} and
\smash{$\mathbb{P}_{p,r}^{\tilde {G}}( \tilde{E}^*_{e_1}
  )=\mathbb{P}_{p,r}( \tilde{E}^*_{e_1} )$}. (These equalities hold
despite the fact that \smash{$\mathbb{P}_{p,r}^{\tilde {G}}$} is
\emph{not} the same distribution as the projection of
$\mathbb{P}_{p,r}$ on $\tilde {G}$.)  Therefore, by Propositions
\ref{fscpropositionforp>0} and \ref{fscstarpropositionforp>0},
\begin{equation}\label{GtildeEe}
  \mathbb{P}_{p,r}^{\tilde {G}}( \tilde{E}_{e_1} )\geq 0.8639
\end{equation}
would imply that $r_{c}(p)\leq r$, and
\begin{equation}\label{GtildeEestar}
  \mathbb{P}_{p,r}^{\tilde {G}}( \tilde{E}^*_{e_1} )\geq 0.8639
\end{equation}
would imply that $r_c(p)\geq 1-r$.  Below we shall describe a method
which tests whether (\ref{GtildeEe}) or (\ref{GtildeEestar}) holds,
simultaneously for all values of $r\in [0,1]$.

We construct the DaC model on $\tilde {G}$ with parameters $p$ and an
arbitrary $r\in [0,1]$ as follows. Fix an arbitrary deterministic
enumeration $v_1,v_2,\ldots ,v_{|\tilde {\mathcal{V}}|}$ of the vertex
set $\tilde {\mathcal{V}}$, and for $V\subset \tilde {\mathcal{V}}$,
let $\min (V)$ denote the vertex in $V$ of the smallest index.  For
all $r\in [0,1]$, we define the function
\begin{equation*}
  \begin{array}{lccc}
    \Psi _r\: : \: & \Omega _{\tilde {\mathcal{E}}} \times [0,1]^{\tilde {\mathcal{V}}} & \rightarrow & \Omega_{\tilde {\mathcal{E}}} \times \Omega_{\tilde {\mathcal{V}}},\\
    &(\eta,U )& \mapsto & (\eta,\xi _r),\\
  \end{array}
\end{equation*}
where
\begin{displaymath}
  \xi _r(v)=
  \left\{ \begin{array}{ll}
      1 &\textrm{if } U(\min(C_v(\eta )))<r,\\
      0 &\textrm{if } U(\min(C_v(\eta )))\geq r.\\
    \end{array} \right.
\end{displaymath}
Now, if $\mathbb{U}$ denotes uniform distribution on the interval
$[0,1]$ and $(\eta ,U)\in \Omega _{\tilde {\mathcal{E}}} \times
[0,1]^{\tilde {\mathcal{V}}}$ is a random configuration with
distribution $\nu_p^{\tilde {\mathcal{E}}} \otimes \mathbb{U}^{\tilde
  {\mathcal{V}}}$, then it is not difficult to see that $(\eta,\xi
_r)=\Psi _r((\eta ,U))$ is a random configuration with distribution
$\mathbb{P} _{p,r}^{\tilde {G}}$.

We are interested in the following question: for what values of $r$
does $(\eta ,\xi _r)\in \tilde{E}_{e_1} $ (respectively, $(\eta ,\xi
_r)\in \tilde{E}^*_{e_1} $) hold?  The first step is to look at the
edges in $\eta $ in $\tilde R_{e_1}\setminus R_{e_1}$ to see if there
is a bond cluster which connects $R_{e_1}$ and the boundary of $\tilde
R_{e_1}$. If no such connection is found, it is easy to see that there
exists a \textit{threshold value} $r_1=r_1(\eta ,U)\in [0,1]$ such
that for all $r\in [0,r_1)$, $(\eta ,\xi _r)\notin \tilde{E}_{e_1}$,
and for all $r\in (r_1,1]$, we have that $(\eta ,\xi _r)\in
\tilde{E}_{e_1}$.  Indeed, the color configurations are coupled in
such a way that if $r^{\prime }\geq r$ and $(\eta ,\xi _r)\in
\tilde{E}_{e_1}$ then $(\eta ,\xi _{r^{\prime }})\in \tilde{E}_{e_1}$,
since all vertices that are black in $\xi _r$ are black in $\xi
_{r^{\prime }}$ as well.  A similar argument shows that in case of
$\eta \notin F_{e_1}$, there exists $r_1^*=r_1^*(\eta ,U)\in [0,1]$
such that $(\eta ,\xi _r)\notin \tilde{E}^*_{e_1}$ for all $r\in
[0,r_1^*)$, whereas $(\eta ,\xi _r)\in \tilde{E}^*_{e_1}$ for all
$r\in (r_1^*,1]$.  Otherwise, \emph{i.e.}, if there is a connection in
$\eta $ between $R_{e_1}$ and the boundary of $\tilde R_{e_1}$, we
know that neither of $\tilde{E}_{e_1}$ or $\tilde{E}^*_{e_1}$ has
occurred. Hence, in that case, we define $r_1=r_1^*=1$, which
preserves the above ``threshold value'' properties as
$(r_1,1]=(r_1^*,1]=\varnothing $.

Now, if we want a confidence interval with confidence level
$1-\varepsilon $ where $\varepsilon >0$ is fixed, we choose positive
integers $m$ and $n$ in such a way that the probability of having at
least $m$ successes among $n$ Bernoulli experiments with success
probability $0.8639$ each is smaller than (but close to) $\varepsilon
/2$. For instance, for a $99.9999\%$ confidence interval, we can
choose $n=400$ and $m=373$.  By repeating the above experiment $n$
times, each time with random variables that are independent of all the
previously used ones, we obtain threshold values $r_1,r_2,\ldots
,r_{n}$ and $r_1^*,r_2^*,\ldots ,r_{n}^*$. Then we sort them so that
$\tilde{r}_{1} \leq \tilde{r}_{2} \leq ... \leq \tilde{r}_{n}$, and
$\tilde{r}_{1}^* \leq \tilde{r}_{2}^* \leq ... \leq \tilde{r}_{n}^*$.
\begin{Proposition}\label{confintprop}
  Each of the inequalities $r_c(p)\leq \tilde{r}_{m}$ and
  $1-\tilde{r}^*_{m}\leq r_c(p)$ occurs with probability at least
  $1-\varepsilon /2$, hence $[1-\tilde{r}^*_{m},\tilde{r}_{m}]$ is a
  confidence interval for $r_c(p)$ of confidence level $1-\varepsilon
  $.
\end{Proposition}

Before turning to the proof, we remark that the above confidence
interval does not necessarily provide meaningful information.  In
fact, with very small ($<\varepsilon $) probability,
$\tilde{r}_{m}<1-\tilde{r}^*_{m}$ can occur. Otherwise, for
unreasonable choices of $s$ and $\ell $, taking a too small $\ell $ in
particular, it could happen that there is a connection in the bond
configuration between $R_{e_1}$ and the boundary of $\tilde R_{e_1}$
in at least $n-m+1$ experiments out of the $n$, in which case
$[1-\tilde{r}^*_{m},\tilde{r}_{m}]=[0,1]$ indeed contains $r_c(p)$ but
gives no new information.

However, the real difficulty is that although a confidence interval
with an arbitrarily high confidence level may be obtained with the
above algorithm, we do not know in advance how \emph{wide} the
confidence interval is. The width of the interval depends on $s$ and
$\ell $, and it is a difficult problem to find good parameter values.
A way to make the confidence interval narrower is to decrease the
value of $m$, but that comes at the price of having a lower confidence
level.

The choices we made for the parameters $s$ and $\ell $ in our
simulations, together with some intuitive reasoning advocating these
choices, are given in the Appendix.  \vspace{0.4cm}

\noindent {\bf Proof of Proposition \ref{confintprop}.}
Let $\mathbb{S}$ be the probability measure on the sample space
$[0,1]^{2n}$ which corresponds to the above experiment, where a
realization
$(\tilde{r}_1,\tilde{r}^*_1,\tilde{r}_2,\tilde{r}^*_2,\ldots
,\tilde{r}_{n},\tilde{r}^*_{n})$ contains the (already ordered)
threshold values. Let $\mathbb{B}_{0.8639}$ denote the binomial
distribution with parameters $n$ and $0.8639$, and $\mathbb{B}_{a(r)}$
the binomial distribution with parameters $n$ and
$a(r)=\mathbb{P}_{p,r}^{\tilde {G}}(\tilde{E}_{e_1})$.

For $r\in [0,1]$, let $N_r$ denote the number of trials among the $n$
such that $\tilde{E}_{e_1}$ occurs at level $r$. Note that $N_r$ has
distribution $\mathbb{B}_{a(r)}$.  Since $a(r)\geq 0.8639$ implies
$r\geq r_c(p)$ (see inequality (\ref{GtildeEe})), we have that
$r<r_c(p)$ implies $a(r)< 0.8639$.  Therefore, for all $r<r_c(p)$,
$\mathbb{B}_{a(r)}$ is stochastically dominated by
$\mathbb{B}_{0.8639}$.  This implies that for all $r<r_c(p)$, we have
that
\begin{eqnarray*}
  \mathbb{S}(\tilde{r}_{m}<r ) & \leq  & \mathbb{S}( N_r\geq m ) \\
  & = & \mathbb{B}_{a(r)}(\{ m,m+1,\ldots ,n \}) \\
  & \leq & \mathbb{B}_{0.8639}(\{ m,m+1,\ldots ,n\})\\
  & \leq & \varepsilon /2,
\end{eqnarray*}
by the definition of $m$ and $n$.

Hence, for all $\delta >0$, we have that $ \mathbb{S}(\tilde{r}_{m}<
r_c(p)-\delta )\leq \varepsilon /2, $ which easily implies that $
\mathbb{S}(\tilde{r}_{m}< r_c(p))\leq \varepsilon /2.  $ We also have
$ \mathbb{S}(\tilde{r}^*_{m}<r^*_c(p) )\leq \varepsilon /2 $ by a
completely analogous computation, which implies by equation
(\ref{rcplusrcstaris1}) that $ \mathbb{S}(1-\tilde{r}^*_{m}> r_c(p)
)\leq \varepsilon /2.  $ Therefore,
\begin{displaymath}
  \mathbb{S}(1-\tilde{r}^*_{m}\leq r_c(p)\leq \tilde{r}_{m}) \geq 1-\varepsilon ,
\end{displaymath}
which is exactly what we wanted to prove.  \qed

\section{Results of the simulations}\label{results-section}

\begin{figure}[ht]
  \begin{center}
    \includegraphics[width=.48\hsize]{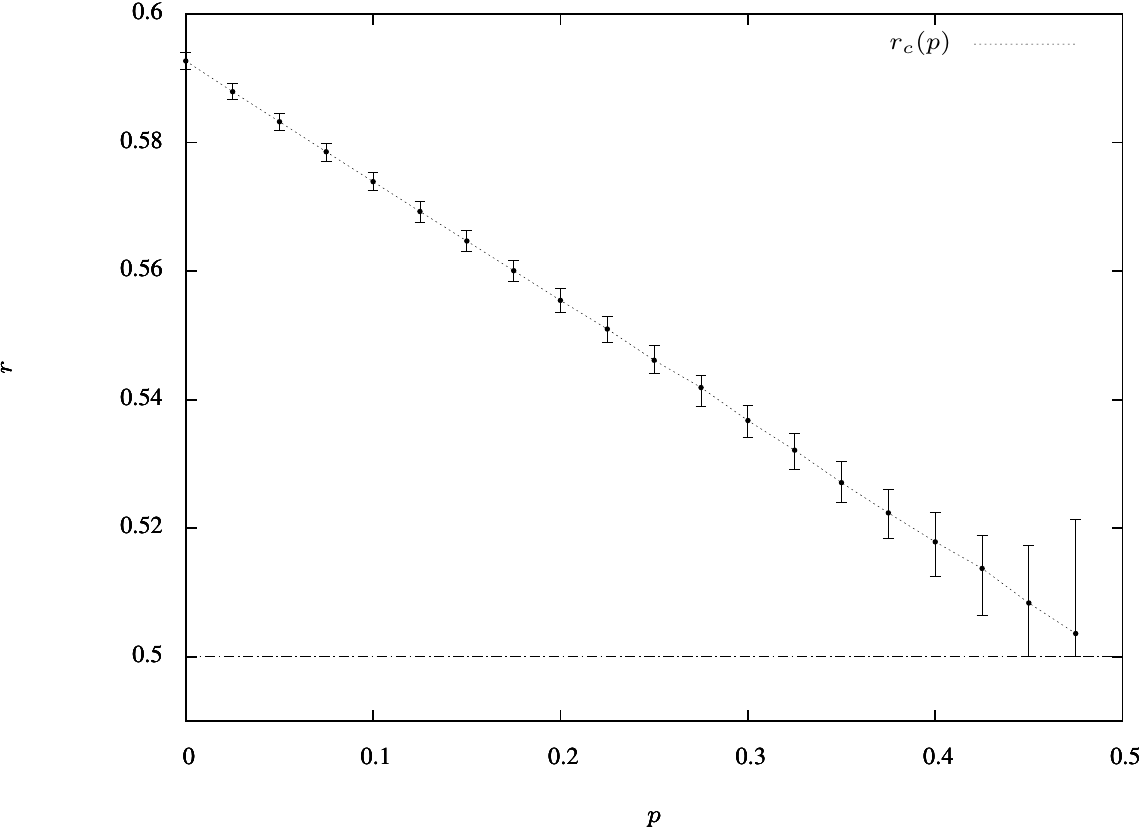}
    \includegraphics[width=.48\hsize]{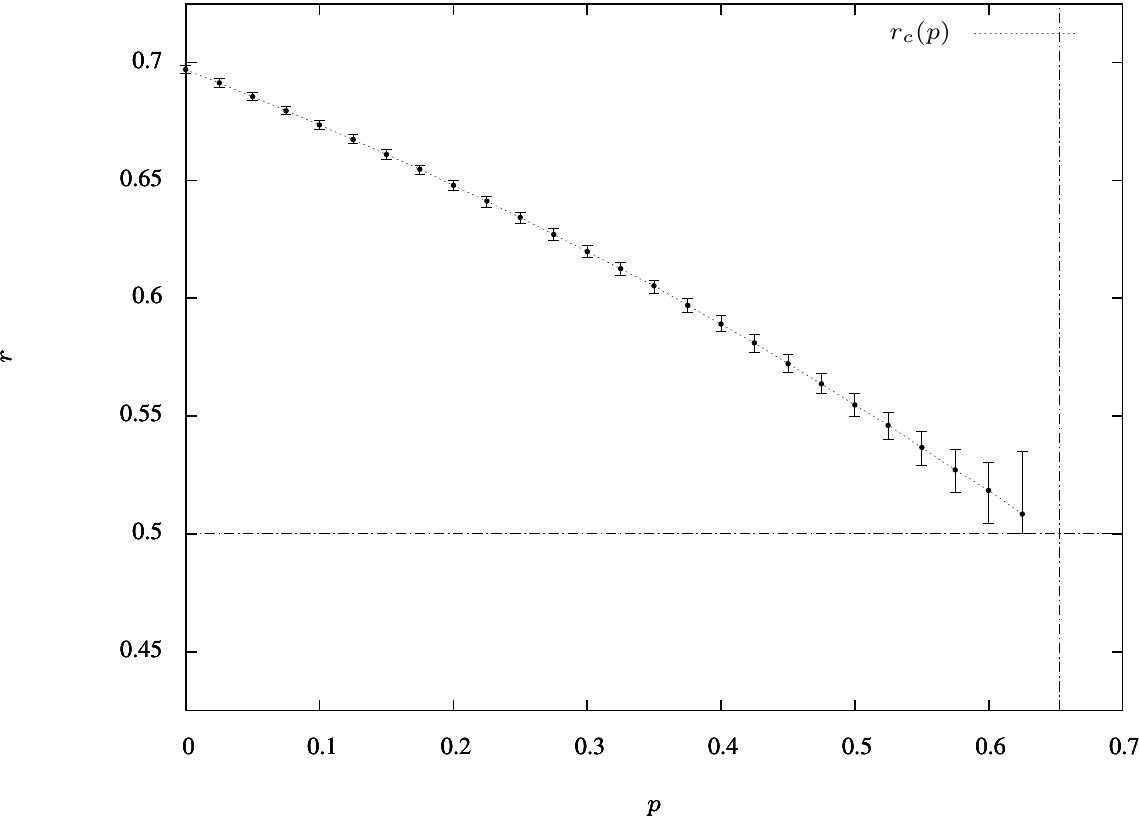}
    \caption{Simulation results for different values of
      $p<p_c^{\mathcal{L}}$ (left: on the square lattice; right: on
      the hexagonal lattice). The dashed line was obtained via a
      non-rigorous correction method.}
    \label{rigorousresults}
  \end{center}
\end{figure}

We implemented the method described in the previous section in a
computer program, and the results for parameter values $\varepsilon
=10^{-6}$, $n=400$, $m=373$ are given below.\footnote{These results
  --- without the description of the method --- have been included in
  \cite{BBT} as well.}  We stress again that although the method in
Section \ref{confidence-interval-section} that determines a confidence
interval for $r_c^{\mathcal{L}}(p)$ is mathematically rigorous, the
results below are obtained by using the random number generator
\cite{mt}, therefore their correctness depends on ``how random'' the
generated numbers are.  The simulations ran on the computers of the
ENS-Lyon, and yielded the confidence intervals represented in
Figure~\ref{rigorousresults}.

Having looked at Figure~\ref{rigorousresults}, we conjecture the
following concerning the behavior of $r_c^{\mathcal{L}}(p)$ as a
function of $p$:

\begin{Conjecture}\label{linearity}
  For $\mathcal{L}\in \{ \mathbb{Z}^2,\mathbb{H}\}$, in the interval
  $p\in [0,p_c^{\mathcal{L}})$, $r_c^{\mathcal{L}}(p)$ is a strictly
  decreasing function of $p$ and $$\lim _{p\to
    p_c^{\mathcal{L}}-}r_c^{\mathcal{L}}(p)= \frac12.$$
\end{Conjecture}

Since it is rigorously known that $r_c^{\mathcal{L}}(0)>1/2$ and
$r_c^{\mathcal{L}}(p)\geq 1/2$ for all $p\in [0,p_c^{\mathcal{L}})$,
Conjecture~\ref{linearity} would imply that $r_c^{\mathcal{L}}(p)>1/2$
for all $p<p_c^{\mathcal{L}}$.  This suggests that the DaC model on
$\mathbb{Z}^2$ or $\mathbb H$ is qualitatively different from the DaC
model on the triangular lattice, where the critical value of $r$ is
$1/2$ for all subcritical $p$ (see Theorem 1.6 in
\cite{BCM}). However, $\lim 
r_c^{\mathcal{L}}(p)=1/2$ would mean that the
difference disappears as $p$ converges to \smash{$p_c^{\mathcal{L}}$}.

The fact that the difference should disappear was conjectured by one
of the authors (VB) and Federico Camia, based on the following
heuristic reasoning.  Near $p=p_c^{\mathcal{L}}$, the structure of the
random graph determined by the bond configuration (whose vertices
correspond to the bond clusters, and there is an edge between two
vertices if the corresponding bond clusters are adjacent in
$\mathcal{L}$) is given by the geometry of ``near-critical percolation
clusters,'' which is expected to be universal for $2$-dimensional
planar graphs. This suggests that the critical $r$ for $p$ close to
its critical value should not depend much on the original underlying
lattice, and we expect the convergence of $r_c^{\mathcal{L}}(p)$ to
$1/2$ to be universal and hold in the case of any $2$-dimensional
lattice.

\bigskip

There is an additional, strange feature appearing in the case of the
square lattice: $r_c(p)$ seems to be close to being an affine function
of $p$ on the interval $[0,1/2)$. This is not at all the same on the
hexagonal lattice, and we have not found any interpretation of this
observation, or of the special role $\mathbb Z^2$ seems to play here.

\begin{Open question}
  Is $r_c^{\mathbb Z^2}(p)$ an affine function of $p$ for $p<1/2$?
\end{Open question}

\paragraph{Acknowledgments.} We thank Federico Camia and Ronald
Meester for interesting and useful discussions and valuable comments
on an earlier version of this paper.  V.T. thanks for the hospitality
of the VU University Amsterdam where much of the work reported here
was derived during an internship provided by ENS Lyon. The research of
A.B. was supported by grants of the Netherlands Organisation for
Scientific Research (NWO) and the Swedish Research Council (VR).
V.B. and V.T. were supported by ANR grant 2010-BLAN-0123-01.

\vspace{0.8cm}
\noindent {\Large {\bf Appendix}}
\vspace{0.5cm}

\noindent The algorithm in Section \ref{confidence-interval-section}
is described for general values of $s$ and $\ell $, and the concrete
values of these parameters will not affect the correctness of the
simulation results.  However, a reasonable choice is important for the
tightness of the bounds obtained and the efficiency of the algorithm,
\emph{i.e.}, the running time of the program.  The heuristic arguments
given here are somewhat arbitrary, and it is quite possible that there
exist other choices that would give at least as good results as ours.

Applying the method described in
Section \ref{confidence-interval-section} requires to simulate a
realization of
the DaC model on the graph $\tilde{G}$,
which is a $2L \times L$ rectangular subset of the square
lattice where
\begin{equation}
  \label{eq:L}
  L=s+2 \ell .
\end{equation}
We will keep this value fixed while we let $\ell $ and $s$
depend on $p$.
Since we want to estimate the critical value for a phase
transition, it is natural to take the largest $L$ possible. After
having performed various trials of our program, we chose
$L=8000$, which was estimated to be the largest value giving a
reasonable time of computation.

Having fixed the size of the graph, we want to choose the parameters
so that the probability of $\tilde {E}_{e_1}$ is as high as possible.
We need to find a balanced value for $\ell $ as small values favor
$E_{e_1}$, but a large $\ell $ might be required to prevent $F_{e_1}$
from happening.  The exponential decay theorem in \cite{ab,menshikov}
for subcritical Bernoulli bond percolation ensures the existence of an
appropriate $\ell $ of moderate size. In our context, we decided that
a good $\ell =\ell (p)$ would be one that ensures
\begin{equation}
  \label{eq:l}
  \mathbb{P}_{p,r}^{\tilde{G}} \left ( F_{e_1} \right ) \approx 0.001.
\end{equation}
We did simulations in order to find an $\ell $ such that (\ref{eq:l}) holds,
then chose $s$ according to equation (\ref{eq:L}). The
values we used in our simulations are summed up in Figure \ref{fig:param}.
\begin{figure}[h]
  \centering
  \begin{tabular}{|c|c|c|c|c|}
\hline
       & \multicolumn{2}{c|}{$\mathbb{Z}^2$}  & \multicolumn{2}{c|}{$\mathbb{H}$}\\
\hline
  $p$    & $s$  & $\ell$ & $s$  & $\ell$ \\
\hline
  $0$    & $7998$ & $1$ &$7998$&$ 1$\\
 \hline
  $0.025$& $7986$ & $7$ &$7986$&$ 7$\\
\hline
  $0.05$ & $7986$ & $8$ &$7984$&$ 8$\\
\hline
  $0.075$& $7982$ & $9$ &$7982$&$ 9$\\
\hline
  $0.1$  & $7980$ & $10$ &$7980$&$ 10$\\
\hline
  $0.125$& $7978$ & $11$ &$7978$&$ 11$\\
\hline
  $0.15$ & $7976$ & $12$ &$7976$&$ 12$\\
\hline
  $0.175$& $7974$ & $13$ &$7974$&$ 13$\\
\hline
  $0.2$  & $7970$ & $15$ &$7970$&$ 15$\\
\hline
  $0.225$& $7964$ & $18$ &$7968$&$ 16$\\
\hline
  $0.25$ & $7962$ & $19$ &$7964$&$ 18$\\
\hline
  $0.275$& $7956$ & $22$ &$7962$&$ 19$\\
\hline
  $0.3$  & $7948$ & $26$ &$7954$&$ 23$\\
\hline
  $0.325$& $7938$ & $31$ &$7952$&$ 24$\\
\hline
  $0.35$ & $7926$ & $37$ &$7946$&$ 27$\\
\hline
  $0.375$& $7904$ & $48$ &$7940$&$ 30$\\
\hline
  $0.4$  & $7876$ & $62$ &$7932$&$ 34$\\
\hline
  $0.425$& $7822$ & $89$ &$7924$&$ 38$\\
\hline
  $0.45$ & $7704$ & $148$ &$7908$&$ 46$\\
\hline
  $0.475$& $7260$ & $370$ &$7896$&$ 52$ \\
\hline
$0.5$ &&& $7876$ & $62$\\
\hline
$0.525$ &&& $7844$ & $ 78$ \\
\hline
$0.55$ &&& $7790$ & $ 105$ \\
\hline
$0.575$ &&& $7710$ & $ 145$ \\
\hline
$0.6$ &&& $7538$ &$231$ \\
\hline
$0.625$ &&& $7002$ &$ 499$ \\
\hline
\end{tabular}
  \caption{Parameters chosen}
  \label{fig:param}
\end{figure}

\end{document}